\documentclass[11pt]{amsart}
 \newtheorem{thm}{Theorem}[section]
 
 \newtheorem{lem}[thm]{Lemma}
 \newtheorem{prop}[thm]{Proposition}
 \theoremstyle{definition}
 
 \theoremstyle{remark}

 \numberwithin{equation}{section}
\usepackage[utf8]{inputenc}   
\usepackage[T1]{fontenc}      
\usepackage{epsfig}
\usepackage{amssymb, graphicx, amsmath, amsfonts}
\usepackage{multicol, hyperref, caption, subcaption, enumerate, longtable, enumitem}
\usepackage{array}
\usepackage{csquotes, rsfso, enumitem, multirow}

\newcolumntype{C}[1]{>{\centering\let\newline\\\arraybackslash\hspace{0pt}}m{#1}}
\usepackage{pifont}

\begin{document}

\title[\tiny Old problem revisited: Which equilateral convex polygons tile?]{Old problem revisited: Which equilateral convex polygons tile the plane? }



\author{Bernhard Klaassen} 

\begin{flushleft}
\center
 Preprint, to appear in Journal for Geometry and Graphics, Vol. 29 (2) 
\end{flushleft}
 
\address{Bernhard Klaassen, Fraunhofer Research Inst. IEG, Bochum, Germany}
\email{bernhard.klaassen@ieg.fraunhofer.de}

\subjclass[2020]{52C20}

\keywords{plane tilings, equilateral tilings, convex tilings}

\date{November 3, 2025}

\begin{abstract}
We present a simplified proof of a forty-year-old result concerning the tiling of the plane with equilateral convex polygons. Our approach is based on a theorem by M. Rao, who used an exhaustive computer search to confirm the completeness of the well-known list of fifteen pentagon types. Assuming the validity of Rao's result, we provide a concise and mainly geometric proof of a tiling theorem originally due to Hirschhorn and Hunt. Finally, a possible connection to quasicrystals is sketched.
\end{abstract}

\maketitle

\section{Introduction}
\label{intro}
The question of which convex equilateral polygons can tile the plane was first thoroughly investigated in 1985 by Hirschhorn and Hunt \cite{Hirsch}. Their approach solved the problem completely, although some parts of their proof relied on computer-aided methods. In 2002, Olga Bagina \cite{Bagina} offered a somewhat shorter proof of the same result, although she also used a small number of symbolic computations (via Maple) to solve certain systems of equations.

Four decades after the original work, further progress has been made, aided in part by computational methods. Notably, Micha\"{e}l Rao \cite{Rao} addressed the broader question of which convex pentagons (with arbitrary edge lengths) tile the plane. He showed, using an exhaustive computer search, that the fifteen known types of tiling convex pentagons form a complete list. Although Rao has not submitted his work to a peer-reviewed journal, Thomas Hales, a highly regarded expert in the field, publicly affirmed the reliability of Rao’s proof in a 2017 blog post \cite{Hales}. If Hales had provided a formal review of Rao’s work for a journal, there would be no reasonable doubt about its correctness. \\
Importantly, Rao's result is independent of the earlier work on equilateral cases by Hirschhorn/Hunt and Bagina.

Assuming that Rao's classification is indeed complete, we can revisit the original problem from a new angle: With the list of fifteen pentagon types accepted as exhaustive, we propose a significantly simpler and mainly geometric approach for characterizing all equilateral convex polygons that tile the plane. The main result of this paper is a geometric proof of the following proposition:

\begin{prop} \label{prop:1}
An equilateral strictly convex polygon P tiles the plane if and only if one of the following holds:
\begin{itemize}[noitemsep, topsep=0pt, label={-}, leftmargin=*]
\item  P is a triangle or quadrangle.
\item  P is a pentagon with two different inner angles that add up to $\pi$ or is similar to a certain pentagon, here called $P7$, later defined in Lemma \ref{lem:3}.
\item  P is a hexagon, and there is a triple of distinct inner angles that add up to $2\pi$ with two of them sharing a common edge.
\end{itemize}
\end{prop}
Proof idea: We will check all 15 known classes\footnote{For completeness, in the appendix all conditions are listed for pentagons and hexagons tiling the plane. During this text, we write ``the polygon \emph{tiles}'' or ``is \emph{tiling}'' for short.} of pentagons and three of hexagons, whether or not they allow convex equilateral solutions and whether all of these valid class representatives fulfill the conditions of Proposition \ref{prop:1}. It will turn out that only very few cases need a closer inspection.\\

Before we proceed with the details of the proof, we should clarify why we only consider strictly convex polygons (and do not include angles equal to $\pi$). It is a convention to exclude polygons with angles equal to $\pi$ if the investigation is restricted to equilateral polygons. This was already the case in the above-mentioned articles \cite{Hirsch} and \cite{Bagina} and even in the dissertation of Reinhardt \cite{Reinhardt}, although he treated polygons with arbitrary side lengths as well. To avoid misunderstandings, we decided to state this explicitly in the proposition.
It should also be noted -- already stated in \cite{Hirsch} -- that any equilateral strictly convex polygon that tiles also admits an edge-to-edge tiling. 

\section{Proofs of lemmas and proposition}

We use the angle notation due to \cite{MMvD} and start the proof section with
\begin{lem}\label{lem:1}
For any convex equilateral pentagon of type $8$ (i.e., $A + B/2 = \pi = C/2 + D$), two inner angles add up to $\pi$.
\end{lem}
\noindent (Proof idea: We construct a symmetric solution of type 8, then check whether two angles add up to $\pi$, then show that no other equilateral convex pentagon is possible for type 8.)
\begin{figure}[!htb]
\centering
\includegraphics[width=1.0\columnwidth]{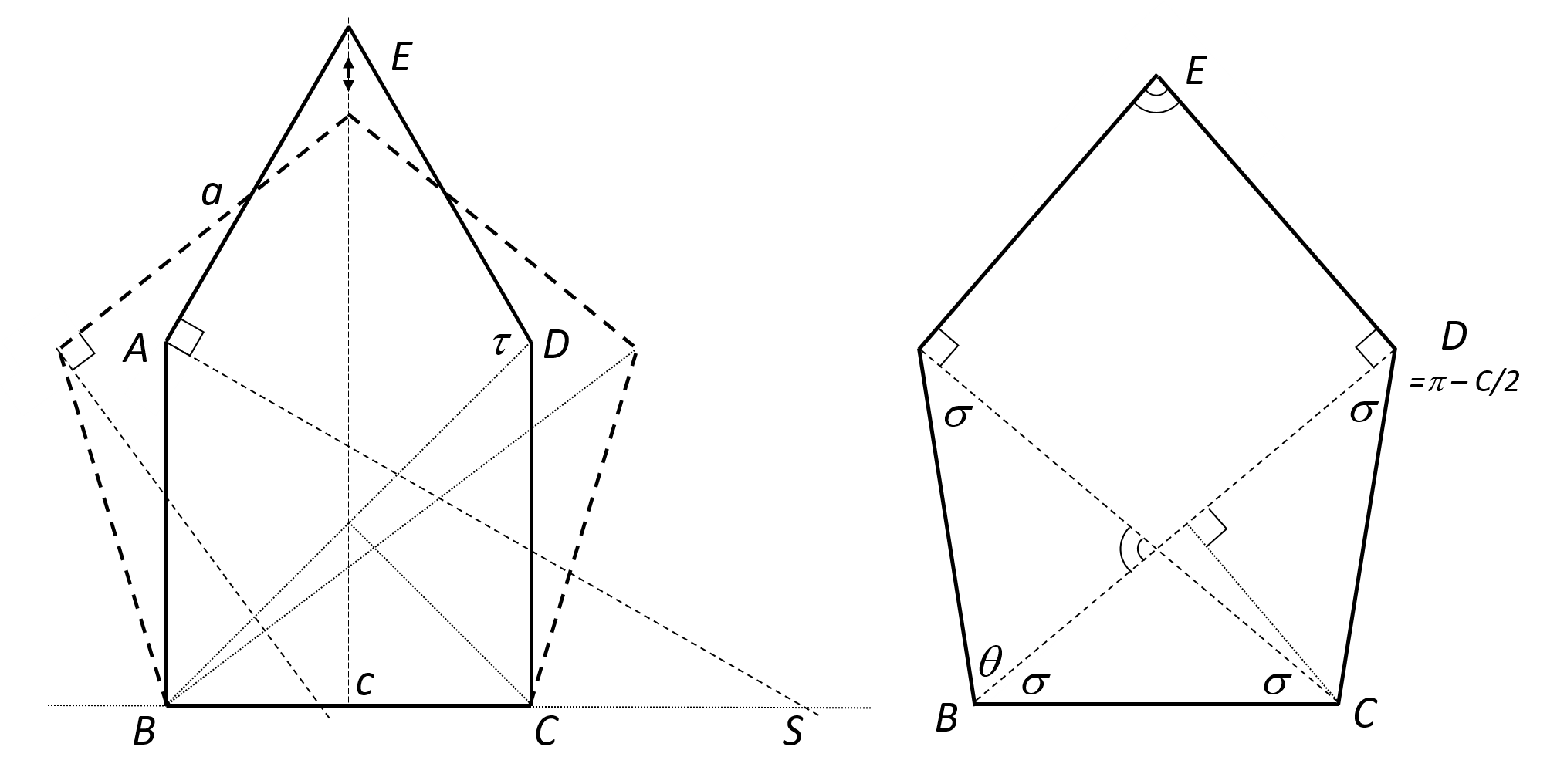}%
\caption{\small Existence of a type $8$ pentagon}%
\label{fig1}%
\end{figure}
\begin{figure}[!htb]
\centering
\includegraphics[width=1.0\columnwidth]{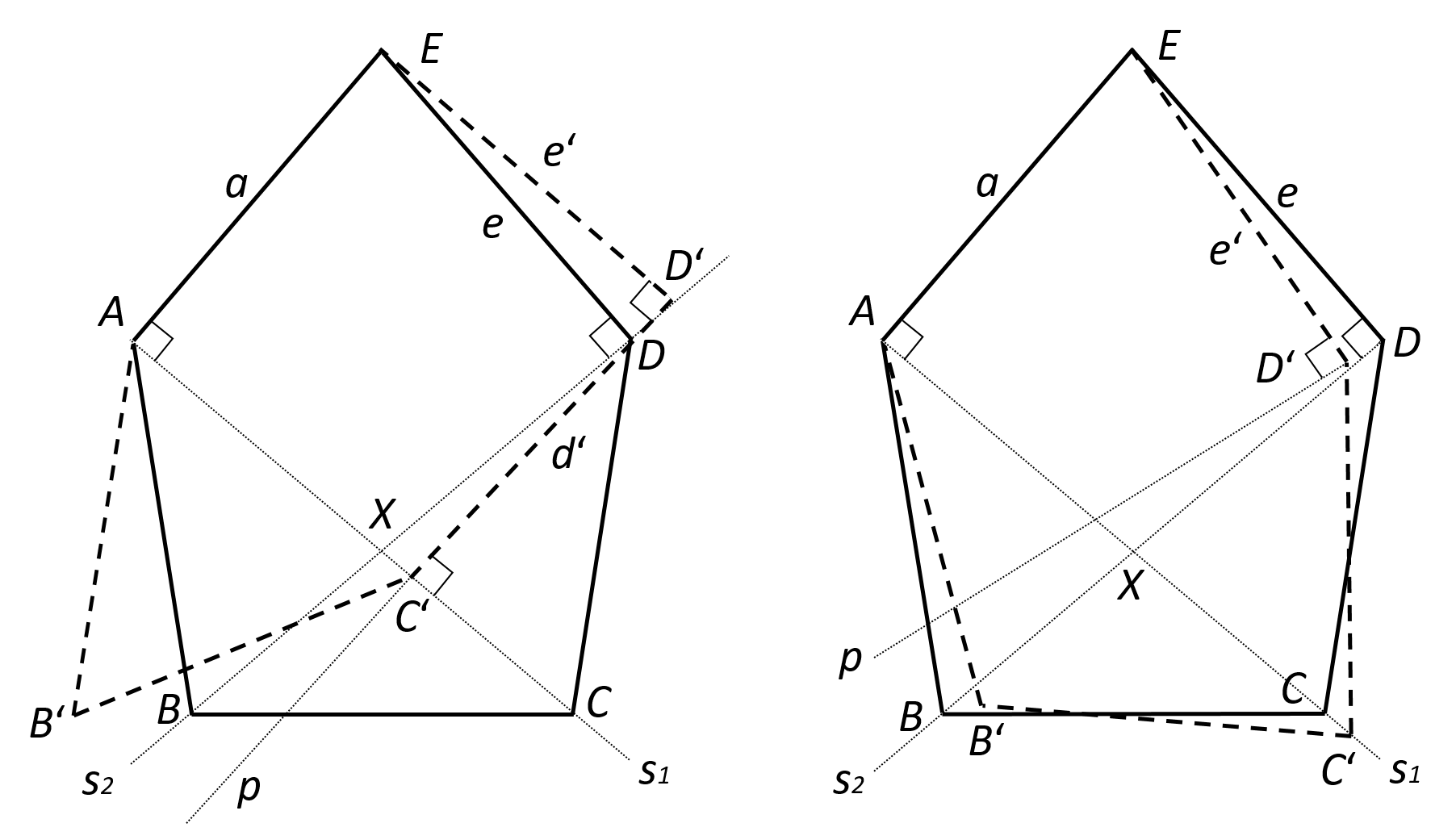}%
\caption{\small Variation of vertex $C$ as $C'$}%
\label{fig2}%
\end{figure}
\begin{proof} 
We set the horizontally fixed edge $c$  (of length $1$), the vertex $E$ shall lie on the perpendicular bisector of $c$ to preserve mirror symmetry (see Fig. \ref{fig1}, left).
Let $\tau$ be the angle formed by the line $\overline{DB}$ and the edge $e$.
Type $8$ condition $C/2+D=\pi$ is equivalent to $\tau=\pi/2$. Consider the following two cases:\\
Case 1: Placing vertex $E$ at distance $1+\sin(60^\circ)$ from $c$ results in $\tau = 105^\circ$.\\
Case 2: Place vertex $E$ such that a regular pentagon is formed with $\tau = 72^\circ$.\\
Hence, there must be a symmetric solution with $\tau=\pi/2$. To make sure that there is only one solution, observe point $S$ at which the perpendicular on edge $a$ crosses the base line. While vertex $E$ is moved downward, edge $a$ is rotated (clockwise around $E$) and translated (downward). Each of these motions is moving $S$ to the left until $C$ is reached. The resulting polygon (see Fig. \ref{fig1}, right) should be called $P8$ and must be strictly convex because both cases 1 and 2 are strictly convex and all angular changes are strictly monotonic under this variation of the vertex $E$. \\
Next, we'll show $B+E=\pi$: Both angles marked with double arcs are equal, since the quadrangle contains two right angles.
The four angles named $\sigma$ are also equal for symmetry reasons.
This implies $B+E=\pi$ since $B=\theta+\sigma$ and $E+\sigma+\theta= \pi$.\\
Now we fix edge $a$ of $P8$ (Fig. \ref{fig2}, left) and check
whether there is a second solution $P'$ with $C'$ different from vertex $C$ in $P8$. ($C'=C$ would imply $P'=P$.)
$C'$ must lie on the line $s_1$ through $A$ and $C$ since $\overline{AC'}$ must be perpendicular to $a$.
If $C'$ is assumed to lie above $C$ on $s_1$, in the section between $C$ and the position of $C'$ shown in the figure, a contradiction is generated since the perpendicular line $p$ on $e'$ crosses the line through $D$ and $B$, called $s_2$ (since $e'$ lies above $e$ during the shift from $C$ to $C'$, because $D>\pi/2$) while $B'$ must lie above $s_2$. 
If $C'$ had been shifted even closer to $A$ than shown in the figure, angle $D'$ would fall below $90^\circ$, and line $p$ would be forced to run outside of $P'$ and could not meet $B'$ as long as $P'$ is convex.\\
On the other hand (Fig. \ref{fig2}, right)
if we assume $C'$ to lie below $C$, then $B'$ must lie below $s_2$ (we have seen in the left part of Fig. \ref{fig2} that the distance $\overline{AX}$ is shorter than a side length). $D'$ now lies above $s_2$ and left of $D$.
Therefore, the perpendicular line $p$ starting at vertex $D'$ to the left stays above $s_2$ and cannot meet $B' \Rightarrow$ contradiction to the assumption that $P'$ exists with $C'$ lower $C$.
Hence, there is no other convex equilateral solution than $P8$ for the type $8$ pentagon.\\
\end{proof} 

\begin{lem}\label{lem:2}
There is no convex equilateral pentagon of type $9$.
\end{lem}
\begin{proof} 
We assume the contrary.  Edge $e$ (see Fig. \ref{fig3}, left) should be vertically fixed. 

\begin{figure}[!htb]
\centering
\includegraphics[width=1.0\columnwidth]{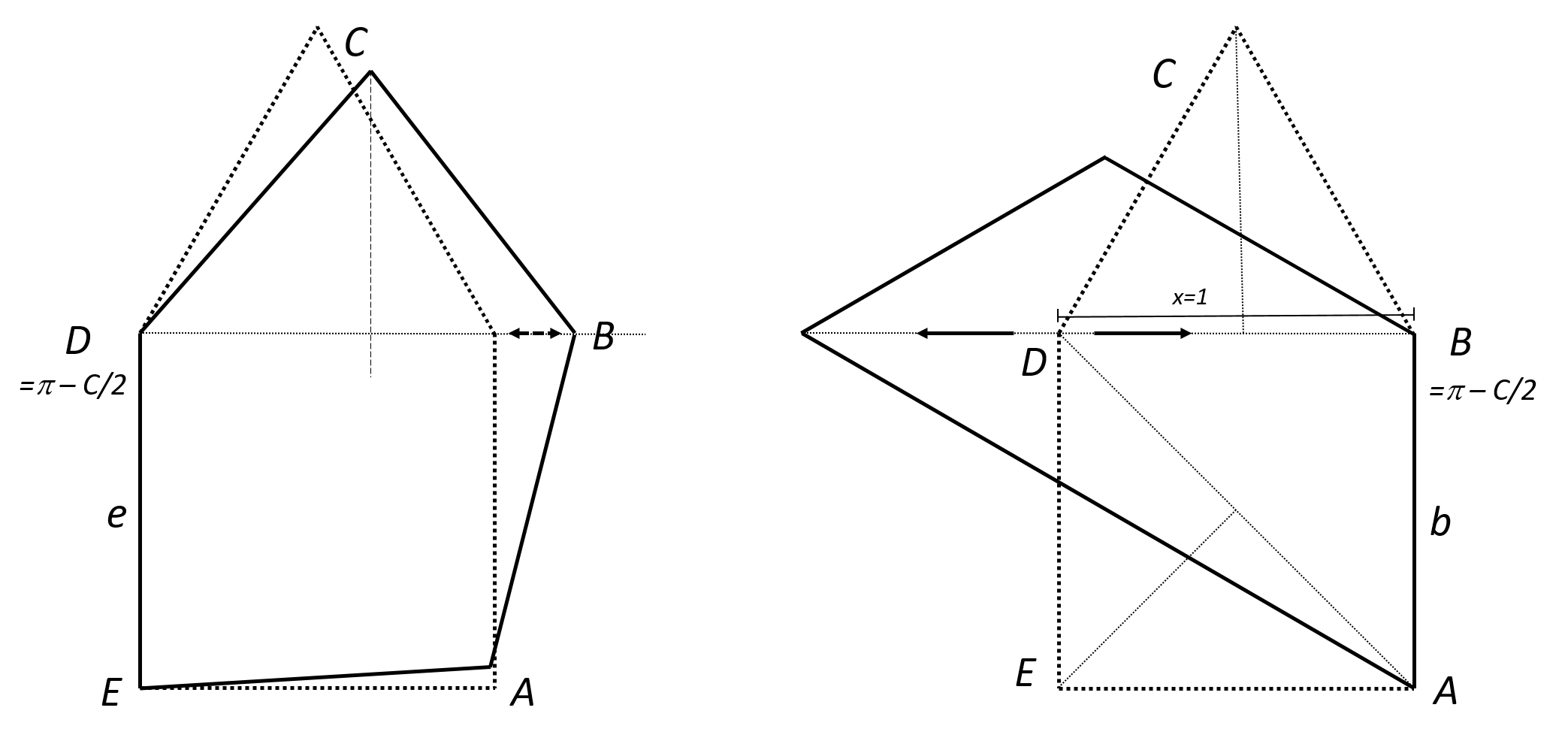}%
\caption{\small Constructing pentagon type $9$ (left) and $7$ (right)}%
\label{fig3}%
\end{figure}
 
\noindent From $D =\pi - C/2$ it follows that $B$ must lie on a line through $D$ perpendicular to $e$;
$\;\;\Rightarrow A$ cannot lie lower than vertex $E$, since $|b|=|e|; \;\;\Rightarrow$ angle $E \leq \pi/2; \Rightarrow$ with $E+B/2=\pi$ it follows that $B/2 \geq \pi/2; 
\;\;\Rightarrow$ due to convexity $B/2= \pi/2; \;\;\Rightarrow E = \pi/2$, satisfied only by the solution with dotted lines; \;\; $\Rightarrow E+B/2<\pi$ which produces a contradiction. 
\end{proof} 

\begin{lem}\label{lem:3}
There exists a unique\footnote{Throughout this paper \emph{unique} means ``unique without regard to similarity''.} convex equilateral representative of type $7$ (called $P7$ in proposition \ref{prop:1})
\end{lem}
\begin{proof} 
Type $7$ needs two conditions: $B+C/2=\pi$ and $A/2+D=\pi$.\\
Fix the edge $b$ (Fig. \ref{fig3}, right) vertically (length = 1). Then vertex $D$ must lie on a horizontal line through $B$ since $B+C/2=\pi$.
Consider $D$ to be placed at a distance $x$ between $0$ and $2\,\cos(30^\circ)$ from $B$. 
For $x=1$ we have $A/2=45^\circ$, $D=150^\circ$, and $A/2+D=195^\circ$ (dotted edges in the figure). At $x=2\,\cos(30^\circ)$ we find  $A/2$ shrinked to $30^\circ$, $D$ to $60^\circ$, the sum to $90^\circ$ (solid edges). It remains to show that the value of $A/2+D$ changes in a strictly monotonic behaviour with $x$, which can easily be seen from the trigonometric functions representing these angles:\\
 $A(x)/2 =0.5\, \operatorname{atan}(x)+ 0.5\, \operatorname{acos}(0.5\sqrt{1+x^2})$ and\\ 
 $D(x) =\operatorname{atan}(1/x)+\operatorname{acos}(0.5\sqrt{1+x^2})+\operatorname{acos}(x/2)$\\
 
 The only term with a positive derivative is $\operatorname{atan}(x)$, all others are strictly monotonic decreasing within the investigated open interval. However, the sum $ \operatorname{atan}(1/x) + 0.5\, \operatorname{atan}(x) $ also results in a monotonic decreasing function since the derivatives of $\operatorname{atan}(x)$ and $\operatorname{atan}(1/x)$ are identical up to a sign. 
 Therefore, it is shown that the unique solution for $x>1$  yields a single equilateral convex pentagon of type $7$, called $P7$ in Proposition \ref{prop:1}.\\
For completeness, we can note here the known numerically calculated angles: $A\!\approx89.26^\circ\!,\, B\!\approx144.56^\circ\!,\, C\!\approx70.88^\circ\!,\, D\!\approx135.37^\circ\!,\, E\!\approx99.93^\circ$. Note that $A<90^\circ$ and $D<150^\circ$, therefore $P7$ is the expected solution with $x>1$ for the equilateral case. (The numerical computation is not part of the proof of this lemma which only says that a unique solution exists, so we can take these values from the literature.)
\end{proof} 

\noindent \emph{Proof of the proposition:}\\
We have to check each type of a tiling polygon. It is well known that Reinhardt \cite{Reinhardt} described all types of strictly convex tiling hexagons and that strictly convex polygons with more than six vertices do not tile. Under the now established assumption that Rao's proof is correct, the list of all candidates is complete. For each type of pentagon or hexagon, one can fill a row in a table (see below) indicating whether it fulfills the criteria of the proposition or whether it can be excluded. We take the angle notation from \cite{MMvD} and \cite{Reinhardt} resp., note that these are equivalent but not coincident with the tables in Wikipedia, but we had to choose a naming convention that is citable and not to be changed in the future. \\
\hspace{1mm}\\
{\small
\begin{tabular}{l l}
  \hline
 Pentagon type & condition (with angle notation due to \cite{MMvD}) \\
  \hline
\textbf{1}  &         $D+E=\pi$  by definition\\ 
\textbf{2} &         $C+E =\pi$  by definition\\
(3)  &        not equilateral by construction (excluded)\\
\textbf{4} &        $A+C =\pi$ since $A=C=\pi/2$  by definition\\ 
\textbf{5} &       $ A+C =\pi$ since $C=2A=2\pi/3$ by definition\\
\hspace{1mm} {\tiny (Note, there is}&\hspace{-5mm}{\tiny an erratum in \cite{MMvD}, fig. 2, type 5: it should read $2\pi/3$ instead of $\pi/2$.)}\\
\textbf{6} &        $C+E =\pi$  by definition\\
\textbf{7} &        Lemma \ref{lem:3}: $P=P7$\\
\textbf{8} &        Lemma \ref{lem:1}:  $B+E=\pi$\\
(9)   &        Lemma \ref{lem:2}: no convex equilateral solution (excluded)\\
(10-15)  &   not equilateral by construction (excluded)\\
  \hline
Hexagon type  & condition (due to Reinhardt \cite{Reinhardt})\\
  \hline
\textbf{1}  &     $A+B+F=2\pi$  by definition\\
\textbf{2}  &      $A+C+F=2\pi$  by definition\\
\textbf{3} &      $B=D=F=2\pi/3 \Rightarrow$ regular hexagon\\
       &   $\Rightarrow$ Any triple of inner angles adds up to $2\pi$.\\
  \hline
\end{tabular}\\
 (Types in bold\footnote{ Note that a type with a number in bold does not guarantee that an equilateral case exists for this type. It just makes sure that if such an equilateral polygon exists, it fulfills the conditions of the proposition. }
 are consistent with the proposition, and the others are excluded.)
 }\\
 
\noindent Finally, we must show that any equilateral strictly convex polygon satisfying the conditions of the proposition always admits a tiling. For triangles, quadrangles, and pentagon $P7$ this is evident. We now assume that such a pentagon exists that does not tile. Then there are only two possible cases: The angles that add up to $\pi$ belong to the same edge, or two edges lie between them. These cases are equivalent to pentagon type $1$ and $2$ resp. \hspace{-0.7mm}in the above table; hence the assumption is wrong. An analogous argument works for hexagons: If two angles of the triple share the same edge, then either the third angle also is sharing an edge with one of the others, or two edges lie between the third angle and one of the two others. Both cases are covered by hexagon type $1$ and $2$ resp. \hspace{-0.7mm}in Reinhardt's list of hexagon types and together with the property that all edges have equal length, it is clear that such polygons tile. This completes the proof.\\

\section{Analysis of the resulting equilateral tilings}
As mentioned in footnote 3, the proposition does not tell us which types of pentagons really exist in equilateral form, it only characterizes a relation between their inner angles. In this section, we will discuss the possible symmetry groups for such tilings. First, we should exclude some further cases for which no equilateral solution exists. 

\subsection{Type 5}
\begin{figure}[!htb]
\centering
\includegraphics[width=1.0\columnwidth]{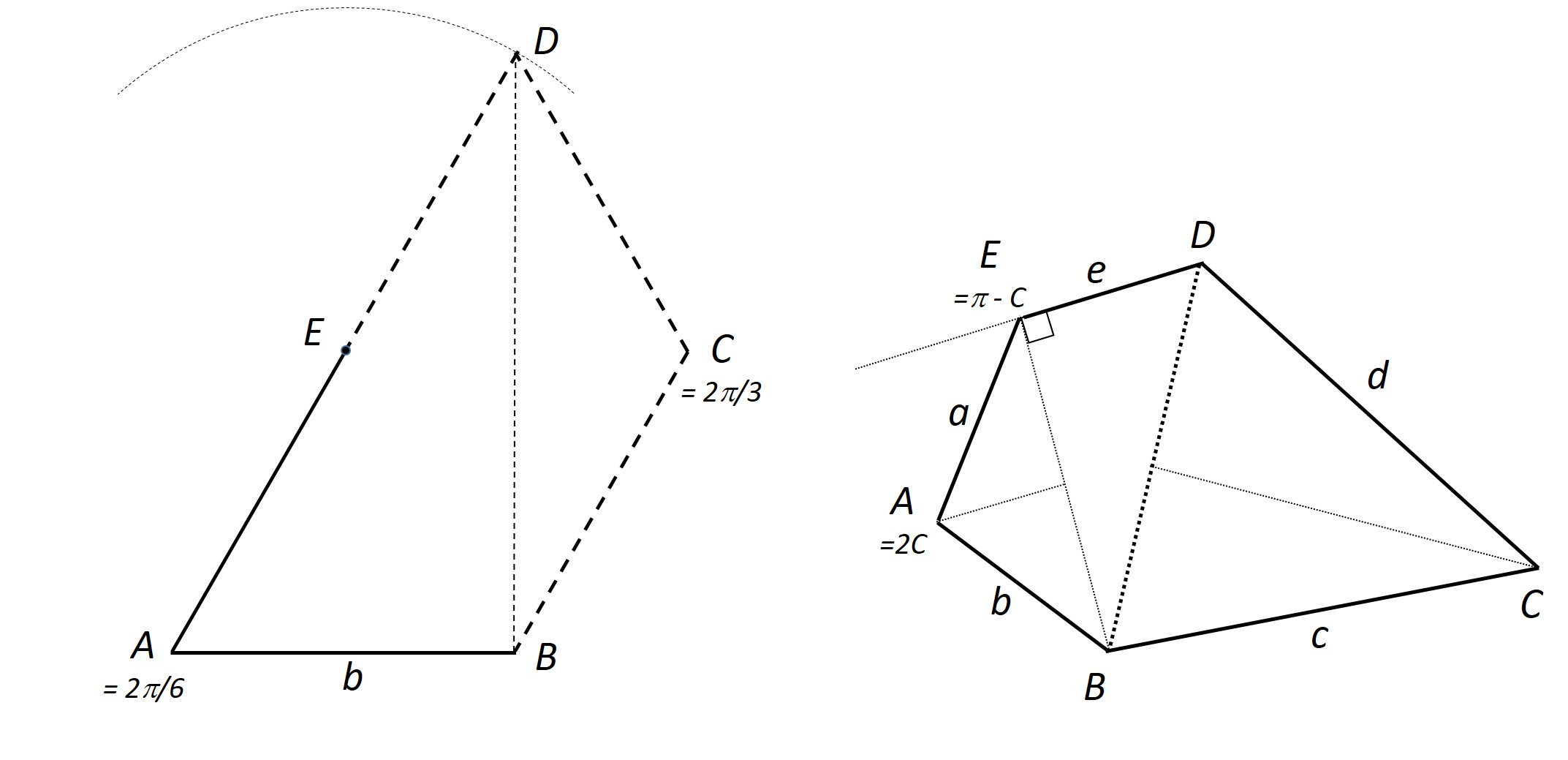}%
\caption{\small Constructing pentagon type $5$ (left) and $6$ (right)}%
\label{fig4}%
\end{figure}

We consider type 5 which has two fixed angles, $A=2\pi/6$ and $C=2\pi/3$ (see fig. \ref{fig4}, left). Observe that with a fixed edge $b$ the position of vertex $E$ is also fixed and $D$ must lie on a circle around $E$. With the predefined size of angle $C$ the distance between $B$ and $D$ is also known. The only possible convex solution is shown in the figure. However, the inner angle $E$ is equal to $\pi$ since all nodes lie on a triangular grid enforced by the two fixed angles $A$ and $C$. Therefore, no strictly convex equilateral solution can exist for type 5.

\subsection{Type 6}
For this type, it is not possible to draw a figure with equal side lengths, so we will demonstrate this fact by a figure showing the conditions for this type (see Fig. \ref{fig4}, right). Observe that edge $e$ lies parallel to the bisecting line of angle $A$ with $|a|=|b|$ and $|c|=|d|$. So, we can express the distance $|BD| = \sqrt{|e|^2 + (2\,|a|\, \sin(C))^2}$ which must be equal to $2\, |d|\, \sin(C/2)$. \newline
We assume that a convex case exists with all side lengths equal to 1. $\sin(C)$ can be expressed by $2 \, \sin(C/2) \cos(C/2)$ and with the substitution $z=\sin^2(C/2)$ we end up with the quadratic equation $z^2 - 3z/4 - 1/16 =0$ which produces only one positive solution $z=(3+\sqrt{13})/8$ and $C\!\approx130.65^\circ$. Hence, $A=2C > \pi$, which is a contradiction to the assumption.

\subsection{Symmetry groups}
Here we list the symmetry groups that occur with equilateral pentagonal tilings. The table is restricted to the periodic edge-to-edge case and we can take the corresponding symmetry groups for each type from the literature \cite{schatt}:\\
\hspace{1mm}\\
{\small
\begin{tabular}{l l}
  \hline
 Pentagon type & possible symmetry groups for periodic equilateral convex tilings \\
  \hline
1  &   \textbf{cmm}, \textbf{cm}, \textbf{pmg} \\ 
2  &   \textbf{pgg} \\
4  &   \textbf{p4} \\
7  &   \textbf{pgg} \\
8  &   \textbf{pgg} \\
  \hline
\end{tabular}\\

 For type 1 several groups are possible. The group with the highest number of different symmetries is \textbf{cmm}, in this sense (not in the sense of group inclusion) this is the maximal symmetry group for type 1.

\section{Nonperiodic equilateral tilings and possible applications}

A naive (but wrong) understanding of the above sections could be that we had found and characterized all possible monohedral tilings with equilateral convex polygons and that all these are periodic. We will demonstrate that such tiles admit several other interesting tilings without a translational symmetry. 
\begin{figure}[!htb]
\centering
\includegraphics[width=0.6\columnwidth]{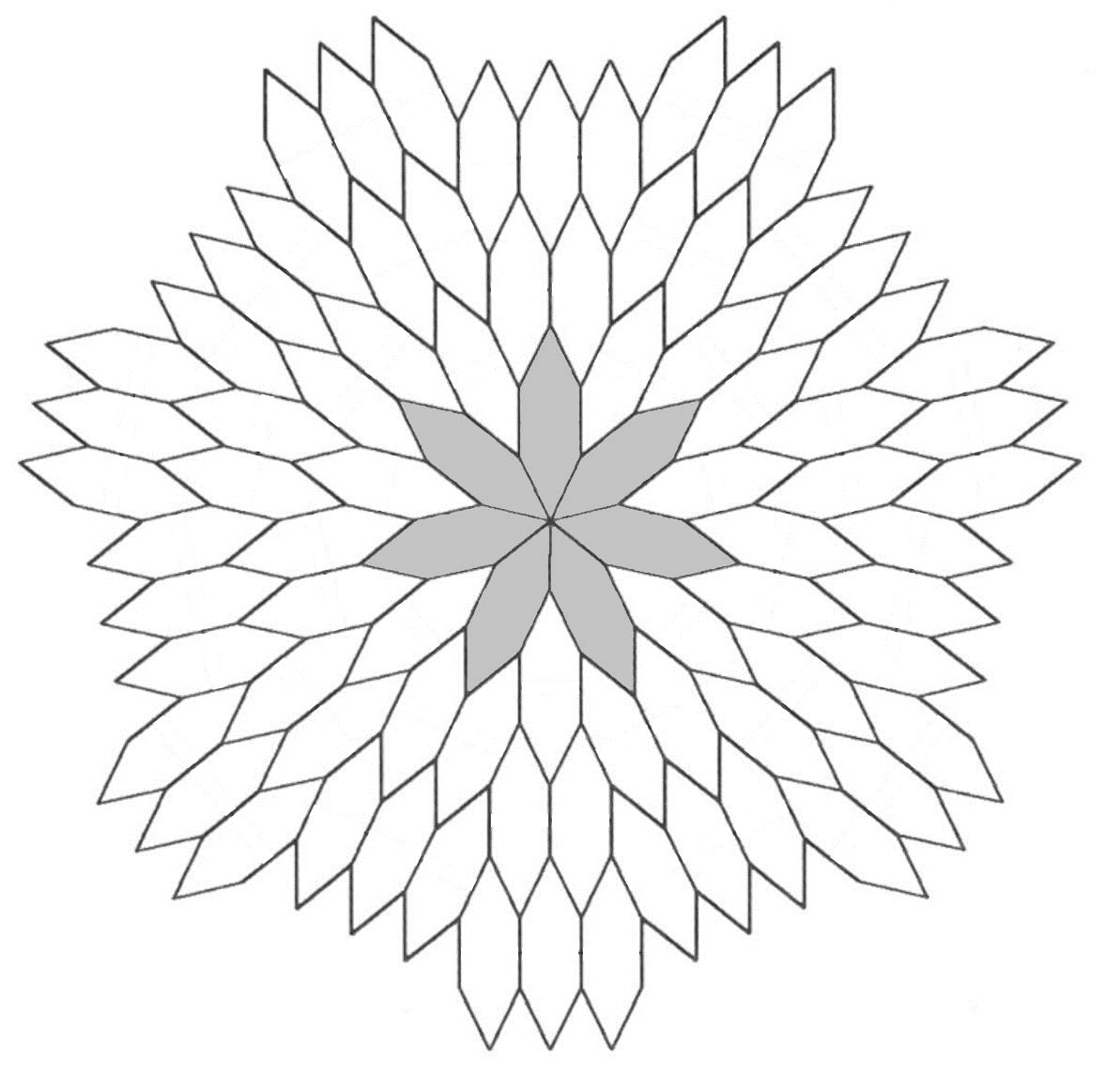}%
\caption{\small Nonperiodic tiling scheme with equilateral hexagons for any rotational symmetry}%
\label{rota7}%
\end{figure}
A classical example with rotational symmetry was already given in the above-mentioned paper by Hirschhorn and Hunt. Moreover, any rotational symmetry of arbitrary order can be represented by an equilateral convex hexagonal tiling, as was demonstrated in \cite{rota}. As an example, Fig. \ref{rota7} shows a tiling with symmetry group $D_7$ in Schoenflies notation.  

Such a nonperiodic tiling should not be misunderstood as aperiodic. Obviously, these monotiles can also generate periodic tilings. Nevertheless, it could be possible that quasicrystals exist with a similar rotational symmetry. Techniques are known to grow crystals starting with a so-called seed. If chemists could create a seed crystal in the form of, e.g., the shaded subset in the center of the Fig. \ref{rota7}, then a periodic crystal structure would be inhibited and only the variant with rotational symmetry could occur. This is, of course, not a guarantee for the existence of quasicrystals but in several cases new materials with unusual rotational symmetries have been found \cite{grow}, \cite{quasi2}. 
\begin{figure}[!htb]
\centering
\includegraphics[width=0.8\columnwidth]{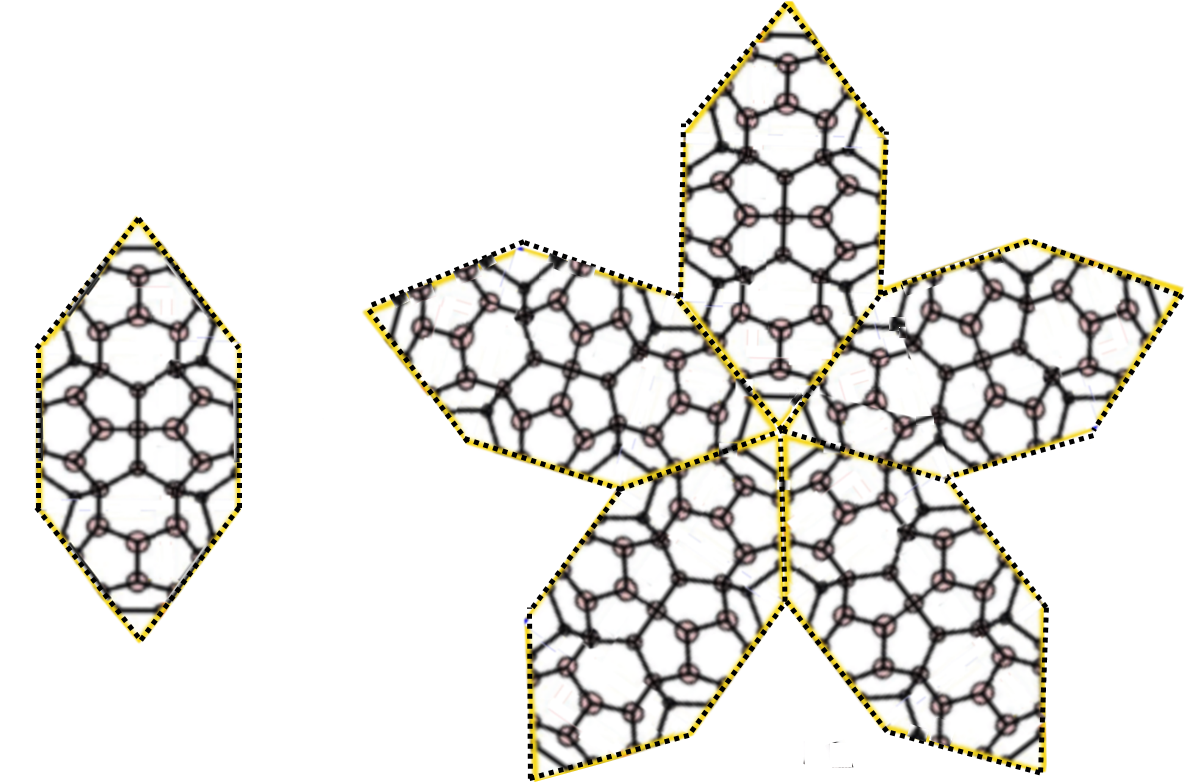}%
\caption{\small Nonperiodic (theoretical) molecular structure composed of equilateral hexagons \cite{quasi3}, marked by dotted lines}%
\label{quasi}%
\end{figure}

Looking further into the molecular structure of \cite{quasi3}, we can even find equilateral hexagonal cells that could (in theory) be composed in the rotational symmetric manner of Fig. \ref{rota7}, here in $D_5$ symmetry. Fig. \ref{quasi} gives an optical impression of this molecular grid. Note that the molecular structure at each of the edges is the same such that any edge-to-edge tiling with this hexagon also works on the molecular level. However, it should be pointed out that the displayed structure does not show an existing quasicrystal. Only the hexagonal cell on the left side of Fig. \ref{quasi} could be found among the examples of \cite{quasi3}.


\section{Alternative approach}
Although - at least to the author's knowledge - there is no serious doubt that Rao's approch was successful and sound, there are  alternative methods to enumerate and categorize tilings.  Meanwhile a well established technique is the Delaney-Dress method \cite{dress}. Here a result of Dress is used which guarantees that two tilings together with their respective symmetry groups are equivalent in structure if and only if the corresponding Delaney symbols for both tilings are isomorphic. This has been utilized in a software-based manner, in particular by Delgado-Friedrichs and Huson \cite{zeller}.

Such methods provide an efficient framework for computer-based analysis of large subsets of tilings under a topological and combinatorical view. However, the method should be enhanced in order to filter the convex cases. Especially the set of pentagonal tilings could be investigated to verify the catalog of 15 convex tiling types in an independent way. So far it is unknown whether this could be done in near future but also other more urgent questions are and will be treated with this framework. It will be interesting to see how this field further develops.

\section{Conclusion}
By accepting that Rao’s result from 2017 on the completeness of the list of 15 tiling pentagons is valid and applying elementary geometric reasoning, we are able to derive a concise characterization of all strictly convex equilateral polygons that admit monohedral tilings of the plane.  In summary, we can give a short and complete description of the resulting polygon types: 
\begin{itemize}[noitemsep, topsep=0pt, label={-}, leftmargin=*]
\item  any equilateral triangle or diamond quadrangle
\item  any strictly convex equilateral pentagon with two different inner angles adding to $\pi$ or pentagon $P7$, defined in Lemma \ref{lem:3}
\item among these pentagons only the types 1, 2, 4, 7 and 8 can occur
\item  any strictly convex equilateral hexagon with three distinct inner angles adding to $2\pi$ while two of them share a common edge.
\end{itemize}

\noindent In addition, some less common nonperiodic tilings were discussed to show the close connection to quasicrystals and to demonstrate the variety of tilings in this special and visually attractive field.


\bibliographystyle{spmpsci}      


%
%

\section*{Appendix: Types of all tiling convex pentagons and hexagons}

Note that these are the conditions for convex pentagons and hexagons with arbitrary side lengths.\\
 \hspace{1cm}\newline
 {\small
\begin{tabular}{l l}
 \hline
  Pentagon & \\
  \hspace{3mm}type & condition(s) with angle notation due to \cite{MMvD} \\
  \hline
1     &         $D+E=\pi$  \\
2     &         $C+E =\pi; \,\,\, a=d$ \\
3     &        $A = C = D = 2\pi/3; \,\,\,a = b, \,\,d = c + e$\\
4     &        $A = C = \pi/2;\,\,\, a = b,\,\, c = d$\\ 
5     &       $C=2A=2\pi/3;\,\,\, a = b, \,\,c = d$\\
\hspace{1mm} {\tiny (Erratum in \cite{MMvD}}&\hspace{-4mm}{\tiny, fig. 2, type 5: it should read $2\pi/3$ instead of $\pi/2$.)}\\6     &        $C + E = \pi, \,\, A = 2C; \,\,\, a = b = e, \,\, c = d$ \\
7     &        $2B + C = 2\pi, \,\, 2D + A = 2\pi; \,\,\, a = b = c = d$\\
8     &        $2A + B = 2\pi, \,\, 2D + C = 2\pi; \,\,\, a = b = c = d$\\
9     &        $2E + B = 2\pi, \,\, 2D + C = 2\pi; \,\,\, a = b = c = d$\\
10   &   $E = \pi/2, \, A + D = 2B - D = \pi, \, 2C + D = 2\pi;  \, a = e = b + d$\\
11  &   $A = \pi/2, \,\, C + E = \pi, \,\, 2B + C = 2\pi; \,\, d = e = 2a + c$\\
12  &   $A = \pi/2, \,\, C + E = \pi, \,\,2B + C = 2\pi; \,\,2a = c + e = d$\\
13  &   $A = C = \pi/2, \,\, 2B = 2E = 2\pi - D; \,\, c = d, \,\, 2c = e$\\
14  &   $D = \pi/2, \,\, 2E + A = 2\pi, \,\, A + C = \pi; \,\, b = c = 2a = 2d$\\
15  &   $A = 60^\circ\!,  B = 135^\circ\!,  C = 105^\circ\!,  D = 90^\circ\!,  E = 150^\circ;  b = d = e, a = 2b$\\
  \hline
Hexagon & \\
\hspace{3mm}type  & conditions (due to Reinhardt \cite{Reinhardt})\\
  \hline
1        &     $A+B+F=2\pi;  \,\, c=f$\\
2       &      $A+C+F=2\pi;  \,\, b=d,\,\, c=f$\\
3       &      $B=D=F=2\pi/3; \,\, a=f, \,\, b=c, \,\, d=e$\\
\end{tabular}
 \hspace{1cm}\newline
 }

\end{document}